\documentclass[11pt]{article}
\usepackage{pslatex}
\usepackage{fancyhdr}
\usepackage{graphicx}
\usepackage{geometry}
\usepackage{amssymb}
\usepackage{amsmath}
\usepackage{graphics}
\usepackage{epsfig}
\usepackage[dvips]{color}
\newtheorem{claim}{\bf \t}[part]

\newtheorem{Definition}{Definition}[part]

\newtheorem{Lemma}{Lemma}[part]
\newtheorem{Proposition}{Proposition}[part]
\newtheorem{Remark}{Remark}[part]
\newtheorem{Theorem}{Theorem}[part]

\numberwithin{Assumption}{section}
\numberwithin{Corollary}{section}
\numberwithin{Definition}{section}
\numberwithin{equation}{section} \numberwithin{Example}{section}
\numberwithin{Lemma}{section} \numberwithin{Proposition}{section}
\numberwithin{Remark}{section} \numberwithin{Theorem}{section}

\def\text#1{{\rm #1}}

\begin{document}
\title{\Large\bf Mean-field backward stochastic differential equations with subdifferrential
operator and its applications  \thanks{The work of Wen Lu is supported
partially by the National Natural Science Foundation of China
(61273128, 11371029) and a Project of Shandong  Province Higher Educational
Science and Technology Program (J13LI06). The work of Yong Ren is
supported by the National Natural Science Foundation of
 China (11371029).  The work of Lanying Hu is
supported by the National Natural Science Foundation of
 China (11201004). }}

\author{\textbf{Wen Lu}$^1$ \footnote{e-mail: llcxw@163.com}\;\  \ \ \textbf{Yong
Ren}$^2$\footnote{Corresponding author. e-mail:
brightry@hotmail.com and renyong@126.com}\;\  \ \  \textbf{Lanying Hu}$^2$\footnote{e-mail: lanyinghu@126.com}
\\
 \small  1. School of Mathematics and Informational Science, Yantai University, Yantai 264005, China   \\
 \small 2. Department of Mathematics, Anhui Normal University, Wuhu
241000, China}\,
\date{}
\maketitle

\begin{abstract}
 In this paper, we deal with a class of mean-field backward stochastic differential equations
  with subdifferrential operator corresponding to a lower semi-continuous convex function. By means
of Yosida approximation, the existence and uniqueness of the solution is established. As an application,
we give a probability interpretation for the viscosity solutions of a class of nonlocal parabolic
 variational inequalities.
\end{abstract}

\vspace{.08in} \noindent \textbf{Keywords:}  Mean-field backward stochastic differential
equation; Subdifferrential operator; Yosida approximation;
McKean-Vlasov equation; Viscosity
solution

\vspace{.08in} \noindent \textbf{\bf MSC} 60H10, 60G40, 60H30

\section{Introduction }

 The general nonlinear case of backward stochastic differential
 equations (BSDEs, in short) was first introduced by Pardoux and Peng \cite{PP}.
 Since then, a lot of works have been devoted to
the study of the theory of BSDEs as well as their applications.
This is due to the connections of BSDEs with mathematical finance, stochastic optimal control as
well as  stochastic games and partial differential equations
(PDEs, in short)  (see e.g.  \cite{EIK}, \cite{HL95}, \cite{HL95-2}
  \cite{Pardoux98}, \cite{Pardoux99}, \cite{PardouxPeng92}, \cite{Peng91}, \cite{Peng93} and so on).

Among the BSDEs, EI Karoui et al. \cite{Karoui1997} introduced a special class of
 reflected BSDEs, which is a
BSDE but the solution $Y$ is forced to stay above a given lower barrier. By means of this kinds of BSDEs, they
 provided a probabilistic formula for the viscosity solution of
an obstacle problem for a class of parabolic PDEs.
In addition, Pardoux and R\c{a}\u{s}canu
[18] proved the existence and uniqueness of the solutions of BSDEs, on a random (possibly
infinite) time interval, involving a subdifferential operator
 (which are also called backward stochastic variational inequalities) in order to give a probabilistic interpretation
for the viscosity solution of some parabolic and elliptic variational inequalities.
Its extension to the probabilistic
interpretation of the viscosity solution of the parabolic
variational inequality (PVI, in short) with a mixed nonlinear
multivalued Neumann-Dirichlet boundary condition was recently given
in Maticiuc and R\u{a}\c{s}canu \cite{MR}.

Mathematical mean-field approaches have been used in many fields,
not only in physics and Chemistry, but also recently in economics,
finance and game theory, see for example, Lasry and Lions
\cite{LasryLions} and the references
therein.

On the other hand, \emph{McKean-Vlasov} stochastic differential equation of the form
\begin{eqnarray*} \label{bsde:1}
  dX_t=b(X_t,\mu_t)dt+dW_t,\quad t\in[0, T],\quad X_0=x,
\end{eqnarray*}
where
$$b(X_t,\mu_t)=\int_{\Omega}b(X_t(\omega), X_t(\omega'))P(d\omega')=E[b(\xi, X_t)]|_{\xi=X_t},$$
$b: R^n\times R^n\rightarrow R$ being a (locally) bounded Borel
measurable function and $\mu(t; \cdot)$ being the probability
distribution of the unknown process $X(t)$, was suggested by Kac
\cite{Kac}  as a stochastic toy model for the Vlasov kinetic
equation of plasma and the study of which was initiated by Mckean
\cite{McKean}. Since then, many authors made contributions on
McKean-Vlasov type SDEs and their applications, see for example,
Ahmed \cite{Ahmed}, Ahmed and Ding \cite{AhmedDing}, Borkar and
Kumar \cite{Borkar}, Chan \cite{Chan}, Crisan and Xiong
\cite{CrisanXiong}, Kotelenez \cite{Kotelenez}, Kotelenez and Kurtz
\cite{KotelenezKurtz} and the references therein.

Recently, Buckdahn et al.
\cite{Buckdahn1} introduced a new kind of BSDEs called as mean-field BSDEs.
Furthermore, Buckdahn et al. \cite{Buckdahn2} deepened the
investigation of mean-field BSDEs in a rather general setting. They
proved the existence and uniqueness as well as a comparison principle
of the solutions  for mean-field BSDEs. Moreover, they established the existence and uniqueness of the viscosity solution for
a class of nonlocal PDEs with the help of the
mean-field BSDE and a McKean-Vlasov forward equation.

Since the works \cite{Buckdahn1} and \cite{Buckdahn2} on the
mean-field BSDEs, there are many efforts devoted to its
generalization. Shi et al. \cite{Shi} introduced and studied
mean-field backward stochastic Volterra integral equations. Xu \cite{Xu} obtained the existence and uniqueness of solutions
for mean-field backward doubly stochastic differential equations, and gave
the probabilistic representation of the solutions for a class of stochastic partial
differential equations by virtue of mean-field BDSDEs. Li
and Luo \cite{LiLuo}  proved the existence and the uniqueness for reflected
mean-field BSDEs. Li \cite{LiJuan} studied reflected mean-filed
BSDEs in a purely probabilistic method, and gave a probabilistic
interpretation of the obstacle problems of the nonlinear and nonlocal PDEs by means of the reflected mean-field BSDEs.

Motivated by the above works, the present paper aims to deal with a class of
mean-field BSDEs with subdifferential operator corresponding to a lower semi-continuous convex function with the following form
\begin{eqnarray}\label{mfsdem:1}\displaystyle
 \left\{ \begin{array}{l@{ }r}-dY_t+\partial\varphi(Y_t)dt\ni E'[f(t,Y_t',Z_t',Y_t,Z_t)]dt-Z_tdW_t,
 \  0\leq t\leq T,\\
\displaystyle Y_T=\xi,
\end{array}\right.
\end{eqnarray}
where $\partial\varphi$ is the subdifferential
operator of a proper, convex and
lower semicontinuous function $\varphi$, $\xi$ is called as the terminal condition.
 More  details refer to Section 2.

The first goal of this paper is to find a triple of adapted processes
$(Y,Z,U)$ in an appropriate space such that mean-field BSDE \eqref{mfsdem:1}  hold
 (see Definition \ref{def:1}). Then, it allow us to establish
 the unique viscosity solution of the following nonlocal parabolic
 variational inequality
 \begin{eqnarray}\label{pvi:1}
 \left\{ \begin{array}{l@{ }r}\displaystyle \frac{\partial u(t,x)}{\partial
 t}+Au(t,x)
 +E[f(t, X_t^{0,x_0}, x, u(t,X_t^{0,x_0}), u(t,x),
 Du(t,x)\cdot E[\sigma(t,X_t^{0,x_0},x)])]\in \partial \varphi(u(t,x)),\\
\displaystyle u(T,x)=E[h(X_T^{0,x_0}, x)], x\in R^n
\end{array}\right.
\end{eqnarray}
with
$$Au(t,x):=\frac{1}{2}tr(E[\sigma(t, X_t^{0,x_0}, x)]E[\sigma(t,
X_t^{0,x_0}, x)]^T D^2u(t,x))
+\langle E[b(t, X_t^{0,x_0}, x)], Du(t,x)\rangle.$$

The paper is organized as follows. In Section 2, we introduce some
notations, basic assumptions and preliminaries. Section 3 is devoted to the proof of the existence
and uniqueness of the solution to the mean-field BSDEs with subdifferential operator by means of Yosida
approximation. In Section 4, we give a probability interpretation for the viscosity solution of a class of nonlocal
 parabolic variational inequalities by means of the
mean-field BSDEs with subdifferential operator.

\section{Notations, preliminaries and basic assumptions}

Let $T>0$ be a fixed deterministic terminal time. Let
$\{W_t\}_{t\geq0}$ be a $d$-dimensional standard Brownian motion
defined on some complete probability space $(\Omega, \mathcal {F},
P)$. We denote by $\mathbb{F}=\{\mathcal {F}_s, 0\leq s\leq T\}$ the
natural filtration generated by $\{W_t\}_{0\leq t\leq T}$ and
augmented by all $P$-null sets, i.e.,
$$\mathcal {F}_s=\sigma\{W_r, r\leq s\}\vee \mathcal {N}_{P}, s\in[0,T],$$
where $\mathcal {N}_{P}$ is the set of all $P$-null subsets. For any
$n\geq1$, $|z|$ denotes the Euclidean norm of $z\in R^n$.

In what follows, we need the following spaces.
 \begin{itemize}
\item  $S^{2}_{\mathbb{F}}(0,T; R)$: the space of $\mathbb{F}$-adapted processes  $Y:
\Omega\times [0,T]\rightarrow  R$ such that
$E\left[ \sup_{t\in[0, T]}|Y_t|^2 \right]<+\infty$;

\item $H^2_{\mathbb{F}}(0,T; R^{n})$: the space of
$\mathbb{F}$-progressively measurable processes $\psi:
\Omega\times[0,T]\rightarrow R^{n} $ such that
 $\|\psi\|^2:= E\int_{0}^T|\psi_t|^2dt<+\infty.$
 \end{itemize}

Now given a measurable function $f: \Omega\times[0,T]\times R\times
R^d\rightarrow R$  which satisfies that $(f(t, y, z))_{t\in[0,T]}$
is $\mathbb{F}$-progressively measurable for all $(y, z)\in R\times
R^d$. We give the following assumption (H0):
\begin{enumerate}
  \item [\rm(i)] There exists some constant $k>0$ such that for all $t\in[0,
T]$, $y, y'\in  R$, $z, z'\in
R^{d}$, it holds that
$$
|f(\omega, t, y,z)-f(\omega, t,
y',z') |  \leq k\Big( |y-y'|+ |z-z'|\Big), \ dP\times dt-{\rm a.s.}
$$
  \item [\rm(ii)] $\displaystyle E\int_0^T|f(s,0,0)|^2ds<+\infty$.
\end{enumerate}
The following result is an immediate consequence of Theorem 1.1 in Pardoux
and R\u{a}\c{s}canu \cite{PardouxRascanu98}.
\begin{Theorem}\label{lemma:2.1}
Let the assumption {\rm(H0)} be satisfied. Then, for the given terminal value $\xi$ satisfying
 $E[|\xi|^2+\varphi(\xi)]<+\infty$,
the BSDE with subdifferential operator
\begin{eqnarray*}
 \left\{ \begin{array}{l@{ }r}-dY_t+\partial\varphi(Y_t)dt\ni f(t,Y_t,Z_t)dt-Z_tdW_t,
 \  0\leq t\leq T,\\
Y_T=\xi,
\end{array}\right.
\end{eqnarray*}
\end{Theorem}
has a unique solution $(Y, Z, U)$ satisfies that
\begin{enumerate}
  \item [\rm(i)] $\displaystyle(Y, Z, U)\in S^{2}_{\mathbb{F}}(0,T; R)\times H^2_{\mathbb{F}}(0,T; R^{d})
\times H^2_{\mathbb{F}}(0,T; R)$.
  \item [\rm(ii)] $\displaystyle  E\int_0^T\varphi(Y_t)dt<+\infty$.
  \item [\rm(iii)] $\displaystyle (Y_t,U_t)\in \partial\varphi$,  $dP\times dt$-{\rm a.e.} {\rm on} $\Omega\times[0,T]$.
  \item [\rm(iv)] $\displaystyle Y_t+\int_t^TU_sds=\xi+\int_t^Tf(s,Y_s,Z_s)ds-\int_t^TZ_sdW_s, 0\leq t \leq T.$
  \end{enumerate}

\subsection{Mean-field BSDEs}

This subsection is devoted to recall some basic results on mean-field BSDEs. The reader is referred
to Buckdahn et al. \cite{Buckdahn1,Buckdahn2} for more details.

Let  $(\bar{\Omega}, \bar{\mathcal{F}},
\bar{P})=(\Omega\times\Omega, \mathcal{F}\otimes\mathcal{F},
P\otimes P )$ be the (non-completed) product of  $(\Omega,
\mathcal{F}, P)$ with itself. We denote the filtration of this
product space by
$\bar{\mathbb{F}}=\{\bar{\mathcal{F}}_t=\mathcal{F}\otimes\mathcal{F}_t,
0\leq t\leq T\}$. A random variable $\xi\in L^{0}(\Omega,
\mathcal{F}, P;  R^{n})$ originally defined on $\Omega$ is extended
canonically to $\bar{\Omega}: \xi'(\omega',\omega)=\xi(\omega'),
(\omega',\omega)\in \bar{\Omega}=\Omega\times\Omega$. For any
$\theta\in L^{1}(\bar{\Omega}, \bar{\mathcal{F}}, \bar{P};  R)$ the
variable $\theta(\cdot  ,\omega):\Omega\rightarrow  R$  belongs to
$L^{1}(\Omega, \mathcal{F}, P;  R)$, $P(d\omega)$-a.s. We denote its
expectation by
$$ E'[\theta(\cdot,\omega)]=\int_{\Omega}\theta(\omega',\omega)P(d\omega').$$
Notice that $ E'[\theta]= E'[\theta(\cdot,\omega)]\in L^{1}(\Omega,
\mathcal{F}, P;  R)$,  and
$$\bar{E}[\theta]\Big(=\int_{\Omega}\theta d\bar{P}
=\int_{\Omega} E'[\theta(\cdot,\omega)]P(d\omega)\Big)= E[
E'[\theta]].$$

Consider the mean-field BSDE
\begin{eqnarray}\label{mfbsde:1}
Y_t=\xi+\int_t^TE'[f(s,Y'_s,Z'_s,Y_s,Z_s)]ds-\int_t^TZ_sdW_s, \
0\leq t\leq T,
\end{eqnarray}
where the driver $f: \bar{\Omega}\times  R\times R^{d}\times  R\times  R^{d}\rightarrow  R$ is
$\bar{\mathbb{F}}$-progressively measurable, and satisfies the
following assumptions.
\begin{enumerate}
  \item [\rm(H1)] There exists a constant  $k>0$ such that for all $t\in[0,
T]$, $y,y',\xi,\xi'\in  R$, $z,z',\eta,\eta'\in
R^{d}$,  we have
\begin{eqnarray*}
&& |f(\omega',\omega, t, y,z,\xi,\eta)-f(\omega',\omega, t,
y',z',\xi',\eta') | \nonumber\\&& \leq k\Big( |y-y'| +
 |\xi-\xi'|+ |z-z'|+
 |\eta-\eta'| \Big), \ \ dP\times dt-{\rm a.s.},
\end{eqnarray*}
  \item [\rm(H2)] $\displaystyle \bar{E}\int_0^T|f(s,0,0,0,0)|^2ds<+\infty$.
\end{enumerate}
\begin{Remark}\rm We emphasize that, due to our notations, the driving coefficient $f$ of
\eqref{mfbsde:1} has to be interpreted as follows
\begin{eqnarray*}
 E'[f(s,Y_s',Z_s',Y_s,Z_s)](\omega)&=&E'[f(s,Y_s',Z_s',Y_s(\omega),Z_s(\omega))]
\\&=&\int_{\Omega}f(s,Y_s'(\omega'),Z_s'(\omega'),Y_s(\omega),Z_s(\omega))P(d\omega').
\end{eqnarray*}
\end{Remark}

The following well-known result is from Buckdahn et al.
\cite{Buckdahn2}.
\begin{Lemma}
Under the assumptions {\rm(H1)} and {\rm(H2)}, for any given $\xi\in
 L^2(\Omega, \mathcal {F}_T, P;  R)$, the mean-field BSDE
\eqref{mfbsde:1} has a unique solution
$(Y,Z)\in S^{2}_{\mathbb{F}}(0,T; R)\times H^2_{\mathbb{F}}(0,T; R^{d})$.
\end{Lemma}

In Buckdahn et al. \cite{Buckdahn2}, the authors also
presented the following comparison result.
\begin{Lemma}
Let $f_i=f_i(\omega,\omega',t,y',z',y,z), i=1,2$ be two drivers
satisfying the assumptions {\rm(H1)} and {\rm(H2)}.  Moreover, we suppose that
\begin{enumerate}
  \item [\rm(i)] one of the two coefficients is independent of $z'$;
  \item [\rm(ii)]  one of the two coefficients is nondecreasing in $y'$.
  \end{enumerate}

Let $\xi_1, \xi_2\in  L^2(\Omega, \mathcal {F}_T, P;  R)$ and
denote by $(Y^1,Z^1)$, $(Y^2,Z^2)$ the solutions of mean-field BSDE
\eqref{mfbsde:1} with data $(\xi_1, f_1)$ and $(\xi_2,
f_2)$, respectively. Then if $\xi_1\geq \xi_2$, $P$-a.s., and
$f_1\geq f_2$, $\bar{P}$-a.s., it holds that also $Y_1\geq Y_2$,
$t\in[0,T]$, $P$-a.s.
\end{Lemma}

\subsection{Mean-field BSDEs  with subdifferential operator}

In this subsection, we introduce some preliminaries of mean-field BSDEs with
subdifferential operator.

Consider the mean-field
BSDE as the form
\begin{eqnarray}\label{mfbsde:2}
 \left\{ \begin{array}{l@{ }r}-dY_t+\partial\varphi(Y_t)dt\ni E'[f(t,Y_t',Z_t',Y_t,Z_t)]dt-Z_tdW_t,
 \  0\leq t\leq T,\\
Y_T=\xi,
\end{array}\right.
\end{eqnarray}
where $\xi$ is  the terminal value and satisfies that
\begin{enumerate}
  \item [\rm(H3)] $\displaystyle E[|\xi|^2+\varphi(\xi)]<+\infty$.
\end{enumerate}
Moreover, $\partial\varphi$ in mean-field
BSDE \eqref{mfbsde:2} is the subdifferential
operator of the function $\varphi: R\rightarrow [0,
+\infty]$ which satisfies the following assumptions:

\begin{enumerate}
  \item [\rm(A1)] $\varphi$ is a proper ($\varphi\not\equiv +\infty$), convex and
lower semicontinuous function,
  \item [\rm(A2)] $\varphi(y)\geq\varphi(0)=0$.
\end{enumerate}

Let us define
\begin{eqnarray*}
&&{\rm Dom}\varphi=\{u\in R: \varphi(u)<+\infty\},\\
&&\partial\varphi(u)=\{u^*\in R: \langle u^*,
v-u\rangle+\varphi(u)\leq\varphi(v), \forall v\in R\},\\
&&{\rm Dom}(\partial\varphi)=\{u\in R:
\partial\varphi(u)\neq\emptyset\},\\
&&(u, u^*)\in \partial\varphi \Longleftrightarrow u\in {\rm
Dom}(\partial\varphi), u^*\in \partial\varphi(u).
\end{eqnarray*}
\begin{Remark}\rm
It is well known that the subdifferential operator $\partial\varphi$ is a
maximal monotone operator, i.e., is maximal in the class of operators
which satisfy the condition
$$ \langle u^*-v^*,
u-v\rangle\geq 0,\; \forall (u, u^*), (v, v^*)\in \partial\varphi.$$
\end{Remark}

We end this section by introduce the definition of the solution for the
mean-filed BSDE \eqref{mfbsde:2}.
\begin{Definition}\label{def:1}\rm
 The triple $(Y, Z, U)$ is called as the solution of mean-filed
BSDE \eqref{mfbsde:2} with subdifferential operator if
\begin{enumerate}
  \item [\rm(i)] $\displaystyle(Y, Z, U)\in S^{2}_{\mathbb{F}}(0,T; R)\times H^2_{\mathbb{F}}(0,T; R^{d})
\times H^2_{\mathbb{F}}(0,T; R)$.
  \item [\rm(ii)] $\displaystyle  E\int_0^T\varphi(Y_t)dt<+\infty$.
  \item [\rm(iii)] $\displaystyle (Y_t,U_t)\in \partial\varphi$,  $dP\times dt$-{\rm a.e.} {\rm on} $\Omega\times[0,T]$.
  \item [\rm(iv)] $\displaystyle Y_t+\int_t^TU_sds=\xi+\int_t^TE'[f(s,Y'_s,Z'_s,Y_s,Z_s)]ds-\int_t^TZ_sdW_s, 0\leq t \leq T.$
  \end{enumerate}
\end{Definition}
\section{Existence and uniqueness of the solution}

This section is devoted to prove the existence and uniqueness
of the solution for \eqref{mfbsde:2}.
Firstly, let us propose the main result of this section.
\begin{Theorem}\label{theorem4:1}
Assume that the assumptions {\rm(H1)--(H3)} hold. Then there exists
a unique solution for the mean-field BSDE \eqref{mfbsde:2}.
\end{Theorem}

We mention that our proof is based on the Yosida approximations. For
this purpose, let's introduce an approximation of the function
$\varphi$ by a convex $C^1$-function $\varphi_\epsilon$,
$\epsilon>0$, defined by
\begin{eqnarray}\label{approxi:2}
\varphi_\epsilon(u)&=&\inf\left\{\frac{1}{2\epsilon}|u-v|^2+\varphi(v):
v\in  R\right\}\nonumber\\&=&\frac{1}{2\epsilon}|u-J_\epsilon
u|^2+\varphi(J_\epsilon u),
\end{eqnarray}
where $J_\epsilon u=(I+ \epsilon\partial\varphi)^{-1}(u)$ is called the \emph{resolvent} of the
monotone operator
of $\partial\varphi$. For reader's convenience, we illustrate some properties of this approximation,
 one can see  Brezis \cite{Brezis1973}  for more details.
\begin{Proposition}\label{prop4:1}
 For all $\epsilon, \delta>0$, $u, v\in
 R$, it holds that
\begin{enumerate}
  \item [\rm(i)] $\varphi_\epsilon$ is a convex function with the
gradient  being a Lipschitz function;
  \item [\rm(ii)] $\displaystyle \varphi_\epsilon(u)\leq \varphi(u)$;
  \item [\rm(iii)] $\displaystyle \nabla
\varphi_\epsilon(u)=\partial\varphi_\epsilon(u)=\frac{u-J_\epsilon
(u)}{\epsilon}\in\partial\varphi(J_\epsilon (u))$;
  \item [\rm(iv)] $\displaystyle |J_\epsilon (u)-J_\epsilon (v)|\leq  |u-v|$;
  \item [\rm(v)] $\displaystyle 0\leq \varphi_\epsilon(u)\leq
\langle\nabla\varphi_\epsilon(u), u\rangle$;
  \item [\rm(vi)] $\displaystyle\langle\nabla\varphi_\epsilon(u)-\nabla\varphi_\delta(v),
u-v\rangle\geq-(\epsilon+\delta)\langle\nabla\varphi_\epsilon(u),
\nabla\varphi_\delta(v)\rangle$.
\end{enumerate}
\end{Proposition}
Since our method is based on the Yosida approximations, let
us consider the following mean-field BSDE
\begin{eqnarray}\label{mfbsde:3}
Y^{\epsilon}_t+\int_t^T\triangledown\varphi_{\epsilon}(Y^{\epsilon}_s)dt=
\xi+\int_t^TE'[f(s,Y^{\epsilon\prime}_s,Z^{\epsilon\prime}_s,Y^{\epsilon}_s,Z^{\epsilon}_s)]dt
-\int_t^TZ^{\epsilon}_tdW_t, 0\leq t \leq T.
\end{eqnarray}
Since
$\nabla\varphi_\varepsilon $ is
Lipschitz continuous, it is known from a recent result of Buchdahn et al. \cite{Buckdahn2}
 that the mean-field BSDE \eqref{mfbsde:3} has a unique solution $
(Y^\varepsilon,Z^\varepsilon)\in S^{2}_{\mathbb{F}}(0,T; R)\times H^2_{\mathbb{F}}(0,T; R^{d})   $.

Setting $$U_t^\varepsilon=\nabla
\varphi_\varepsilon(Y_t^\varepsilon), \ 0\leq t \leq T,$$ we shall prove the convergence of the sequence
$(Y^\varepsilon,U^\varepsilon,Z^\varepsilon)$ to a
process $(Y,U,Z)$, which is the desired solution of the mean-field BSDE \eqref{mfbsde:3}.

Firstly, we establish some properties of the solution of mean-field
BSDE \eqref{mfbsde:3}. In what follows, $C>0$ denotes a constant whose value may
change from line to line.
\begin{Lemma}\label{lemma:4.1}
Assume that the assumptions {\rm(H1)--(H3)} hold. Then there exist two positive constants
$\lambda$ and $C$ such that
\begin{eqnarray}
 E\left[\sup_{t\in [0,T]}e^{\lambda
t}|Y^\epsilon_t|^2+\int_0^Te^{\lambda s} |Z^{\epsilon}_s
|^2ds\right]\leq CM_1,
\end{eqnarray}
 where $M_1:=  E[e^{\lambda T}|\xi|^2]+\bar{E}\int_0^Te^{\lambda s} |f(s,0,0,0,0)
|^2ds.$
\end{Lemma}
\textbf{Proof.}  It\^{o}'s formula yields that
\begin{eqnarray}\label{equ:3}
&&e^{\lambda t}|Y_t^\epsilon|^2+\int_t^Te^{\lambda
s}(\lambda|Y_s^\epsilon|^2+|Z_s^\epsilon|^2)ds+2\int_t^Te^{\lambda
s}\langle\nabla\varphi_\epsilon(Y_s^\epsilon),
Y_s^\epsilon\rangle ds\nonumber\\
&=&e^{\lambda T}|\xi|^2+2\int_t^Te^{\lambda s}\langle Y_s^\epsilon,
E'[f(s,Y_s^{\epsilon\,\prime},
Z_s^{\epsilon\,\prime},Y_s^{\epsilon},Z_s^{\epsilon})]\rangle ds\nonumber\\
&&-2\int_t^Te^{\lambda s}\langle
 Y_s^{\epsilon},
Z_s^{\epsilon}dW_s\rangle.
\end{eqnarray}
By Young's inequality and (H1), we have, for $\gamma>0$
\begin{eqnarray}\label{equ:4}
2\langle Y_s^\epsilon, E'[f(s,Y_s^{\epsilon\,\prime},
Z_s^{\epsilon\,\prime},Y_s^{\epsilon},Z_s^{\epsilon})]\rangle
&\leq&\gamma
|Y_s^\epsilon|^2+\frac{1}{\gamma}|E'[f(s,Y_s^{\epsilon\,\prime},
Z_s^{\epsilon\,\prime},Y_s^{\epsilon},Z_s^{\epsilon})]|^2
\nonumber\\
&\leq&\gamma |Y_s^\epsilon|^2+\frac{8
k^2}{\gamma}\Big[E'|Y_s^{\epsilon\,\prime}|^2+E'|Z_s^{\epsilon\,\prime}|^2\Big]
\nonumber\\&&
+\frac{8
k^2}{\gamma}\Big[|Y_s^{\epsilon}|^2+
|Z_s^{\epsilon}|^2\Big]+\frac{2}{\gamma}E'|f(s,0,0,0,0)|^2.
\end{eqnarray}
Since $\langle\nabla\varphi_\epsilon(y), y\rangle\geq 0$, and hence,
\begin{eqnarray}\label{equ:5}
&&e^{\lambda t}|Y_t^\epsilon|^2+\left(\lambda-\gamma-\frac{8
k^2}{\gamma}\right) \int_t^Te^{\lambda
s}|Y_s^\epsilon|^2ds+ \left(1-\frac{8
k^2}{\gamma}\right)\int_t^Te^{\lambda
s}|Z_s^\epsilon|^2ds \nonumber\\
&\leq&e^{\lambda T}|\xi|^2+\frac{8
k^2}{\gamma}\int_t^Te^{\lambda s}
(E'|Y_s^{\epsilon\,\prime}|^2+E'|Z_s^{\epsilon\,\prime}|^2)ds\nonumber\\
&&+\frac{2}{\gamma}\int_t^Te^{\lambda s}
E'|f(s,0,0,0,0)|^2 ds-2\int_t^Te^{\lambda s}\langle
 Y_s^{\epsilon},
Z_s^{\epsilon}dW_s\rangle.
\end{eqnarray}
Choosing $\gamma=16k^2$ and
$\lambda>\gamma+\frac{16 k^2}{\gamma}$, then there exists a constant
$C>0$, depending on $\lambda, \gamma$ and $k$, such that
\begin{eqnarray}\label{equ:6}
E\left[\int_0^Te^{\lambda s} |Y_s^\epsilon|^2 ds+\int_0^Te^{\lambda s} |Z_s^\epsilon|^2 ds\right]
 \leq C\left[ Ee^{\lambda
T}|\xi|^2+\bar{E}\int_0^Te^{\lambda s}|f(s,0,0,0,0)|^2ds\right].
\end{eqnarray}
On the other hand, combining \eqref{equ:5} and \eqref{equ:6}, we get
\begin{eqnarray}\label{equ:7}
\sup_{t\in[0,T]}e^{\lambda t}|Y_t^\epsilon|^2 &\leq& C\Big[e^{\lambda
T}|\xi|^2+\int_0^Te^{\lambda s} E'|f(s,0,0,0,0)|^2ds \nonumber\\
&&+Ee^{\lambda
T}|\xi|^2+\bar{E}\int_0^Te^{\lambda s}|f(s,0,0,0,0)|^2ds\Big]
 \nonumber\\
&&+2\sup_{t\in[0,T]} \left|\int_t^Te^{\lambda s}\langle
 Y_s^{\epsilon},
Z_s^{\epsilon}dW_s\rangle\right|.
\end{eqnarray}
Thus, from Burkholder-Davis-Gundy's inequality, we have
\begin{eqnarray}\label{equ:8}
E\left[\sup_{t\in[0,T]}e^{\lambda t}|Y_t^\epsilon|^2\right] &\leq& C\left[Ee^{\lambda
T}|\xi|^2+\bar{E}\int_0^Te^{\lambda s}|f(s,0,0,0,0)|^2ds\right] \nonumber\\
&&+\frac{1}{2}E\left[\sup_{t\in[0,T]}e^{\lambda t}|Y_t^\epsilon|^2\right]
+CE\int_0^Te^{\lambda s} |Z_s^\epsilon|^2 ds.
\end{eqnarray}
We then complete the proof by \eqref{equ:6}.

\begin{Lemma}\label{lemma:4.2}
Assume that the assumptions {\rm(H1)--(H3)} hold. Then there exists a constant
$C>0$ such that for all $t\in[0,T]$,

\begin{enumerate}
\item [\rm(i)] $\displaystyle E\int_0^T e^{\lambda s}|\nabla\varphi_\epsilon(Y_s^\epsilon)|^2ds\leq
CM_2$,
\item [\rm(ii)] $\displaystyle E e^{\lambda t}\varphi(J_\epsilon (Y_t^\epsilon))+E\int_t^T
e^{\lambda
s}\varphi(J_\epsilon (Y_s^\epsilon))ds\leq CM_2$,
\item [\rm(iii)] $ \displaystyle E(e^{\lambda t}|Y_t^\epsilon-J_\epsilon Y_t^\epsilon|^2)\leq
C\epsilon M_2$,
\end{enumerate}
where $\displaystyle M_2:=E[e^{\lambda T}\varphi(\xi)+e^{\lambda T}|\xi|^2]+\bar{E}\int_0^Te^{\lambda s} |f(s,0,0,0,0)
|^2ds$.
\end{Lemma}
\textbf{Proof.} The stochastic subdifferential inequality in Pardoux
and R\u{a}\c{s}canu \cite{PardouxRascanu98} gives that
\begin{eqnarray*}
e^{\lambda T}\varphi_{\epsilon}(\xi)  \geq  e^{\lambda
t}\varphi_{\epsilon}(Y^{\epsilon}_t)+\int_t^Te^{\lambda
s}\langle\nabla\varphi_{\epsilon}(Y^{\epsilon}_s),
dY^{\epsilon}_s\rangle
+\int_t^T\varphi_{\epsilon}(Y^{\epsilon}_s)d(e^{\lambda s}),
\end{eqnarray*}
and hence
\begin{eqnarray}\label{lm:15}
&& e^{\lambda
t}\varphi_{\epsilon}(Y^{\epsilon}_t)+\lambda\int_t^Te^{\lambda
s}\varphi_{\epsilon}(Y^{\epsilon}_s)ds+\int_t^Te^{\lambda
s}|\nabla\varphi_{\epsilon}(Y^{\epsilon}_s)|^2ds \nonumber\\&\leq&
e^{\lambda T}\varphi_{\epsilon}(\xi) + \int_t^Te^{\lambda s}\langle
\nabla\varphi_{\epsilon}(Y^{\epsilon}_s),
E'[f(s,Y^{\epsilon\prime}_s,Z^{\epsilon\prime}_s, Y^{\epsilon}_s,Z^{\epsilon}_s)]\rangle
ds\nonumber\\&& -\int_t^Te^{\lambda s}\langle
\nabla\varphi_{\epsilon}(Y^{\epsilon}_s), Z^{\epsilon}_sdW_s\rangle.
\end{eqnarray}
By Young's inequality and (H2), we have
\begin{eqnarray}\label{lm:16}
&&\int_0^Te^{\lambda s}\langle
\nabla\varphi_{\epsilon}(Y^{\epsilon}_s),
E'[f(s,Y^{\epsilon\prime}_s,Z^{\epsilon\prime}_s,Y^{\epsilon}_s,Z^{\epsilon}_s)]\rangle
ds\nonumber\\&\leq& \frac{1}{2}\int_0^Te^{\lambda s}
|\nabla\varphi_{\epsilon}(Y^{\epsilon}_s)|^2ds+\int_0^Te^{\lambda s}
E'|f(s, 0, 0, 0, 0)|^2ds \nonumber\\&& +4k^2\int_0^Te^{\lambda s}(|
Y^{\epsilon}_s|^2+E'|Y^{\epsilon\prime}_s|^2+|Z^{\epsilon}_s |^2
+E'|Z^{\epsilon\prime}_s|^2)ds.
\end{eqnarray}
This together with \eqref{lm:15} yields
\begin{eqnarray}\label{lm:17}
 \frac{1}{2}E\int_0^Te^{\lambda s}
|\nabla\varphi_{\epsilon}(Y^{\epsilon}_s)|^2ds&\leq&e^{\lambda
T}\varphi_{\epsilon}(\xi)+\bar{E}\int_0^Te^{\lambda s}|f(s, 0, 0, 0,
0)|^2ds \nonumber\\&&+4k^2E\int_0^Te^{\lambda s}(|
Y^{\epsilon}_s|^2+E'|Y^{\epsilon\prime}_s|^2+|Z^{\epsilon}_s |^2
+E'|Z^{\epsilon\prime}_s|^2)ds.
\end{eqnarray}
Thus, (i) is hold from Lemma \ref{lemma:4.1} and the fact that
 $\varphi_{\epsilon}(u)\leq\varphi(u)$.

On the other hand,  combining \eqref{lm:15} and \eqref{lm:16},
 we get from Lemma \ref{lemma:4.1} that
 \begin{eqnarray}\label{lm:18}
E e^{\lambda t}\varphi_\epsilon (Y_t^\epsilon)+E\int_0^T
e^{\lambda
s}\varphi_\epsilon (Y_s^\epsilon)ds\leq CM_2.
\end{eqnarray}
Since $\varphi(J_\epsilon (y))\leq \varphi_\epsilon (y)$, it follows that
$$E e^{\lambda t}\varphi(J_\epsilon (Y_t^\epsilon))+E\int_0^T
e^{\lambda
s}\varphi(J_\epsilon (Y_s^\epsilon))ds\leq CM_2.$$
For (iii), since
$$\frac{1}{2\epsilon}e^{\beta t}|Y^{\epsilon}_t-J_{\epsilon}(Y^{\epsilon}_t)|^2
\leq e^{\beta t}\varphi_{\epsilon}(Y^{\epsilon}_t),$$ it follows from \eqref{lm:18} that
$$E\left[e^{\beta t}|Y^{\epsilon}_t-J_{\epsilon}(Y^{\epsilon}_t)|^2\right]
\leq C \epsilon M_2.$$
The proof is complete.

\begin{Lemma}\label{lemma:4.3}
Assume that the assumptions {\rm(H1)--(H3)} hold. Then
\begin{eqnarray*}
E\left[\sup_{t\in[0,T]}e^{\lambda
t}|Y_t^\epsilon-Y_t^\delta|^2\right]+E\int_0^Te^{\lambda
t}(|Y_t^\epsilon-Y_t^\delta|^2+|Z_t^\epsilon-Z_t^\delta|^2)dt \leq
C(\epsilon+\delta)M_2.
 \end{eqnarray*}
\end{Lemma}
\textbf{Proof.} By It\^{o}'s formula, we have
\begin{eqnarray}\label{lm1:18}
&&e^{\lambda t}|Y^\epsilon_t-Y^\delta_t|^2+\lambda\int_t^Te^{\lambda
s}|Y^\epsilon_s-Y^\delta_s|^2ds+\int_t^Te^{\lambda s}
|Z^\epsilon_s-Z^\delta_s|^2ds\nonumber\\&&+2\int_t^Te^{\lambda
s}\langle Y^\epsilon_s-Y^\delta_s,
 \nabla\varphi_{\epsilon}(Y^\epsilon_s)-\nabla\varphi_{\delta}(Y^\delta_s) \rangle
ds\nonumber
\\&=& 2\int_t^Te^{\lambda s}\langle
Y^\epsilon_s-Y^\delta_s,
 E'[f(s,Y^{\epsilon\prime}_s,Z^{\epsilon\prime}_s,Y^\epsilon_s,
Z^\epsilon_s)-f(s,Y^{\delta\prime}_s,Z^{\delta\prime}_s,Y^\delta_s,
Z^\delta_s)] \rangle ds  \nonumber
\\&&
-2\int_t^Te^{\lambda s}\langle Y^\epsilon_s-Y^\delta_s,
(Z^\epsilon_s-Z^\delta_s)dW_s\rangle.
\end{eqnarray}
Using Young's inequality and (H2),  we get for $\gamma>0$
\begin{eqnarray}\label{lm1:19}
&& 2\int_t^Te^{\lambda s}\langle Y^\epsilon_s-Y^\delta_s,
 E'[f(s,Y^{\epsilon\prime}_s,Z^{\epsilon\prime}_s,Y^\epsilon_s,
Z^\epsilon_s)-f(s,Y^{\delta\prime}_s,Z^{\delta\prime}_s,Y^\delta_s,
Z^\delta_s)] \rangle ds  \nonumber
\\&\leq&
\gamma\int_t^Te^{\lambda
s}|Y^\epsilon_s-Y^\delta_s|^2ds+\frac{4k^2}{\gamma}\int_t^Te^{\lambda
s}[|Y^\epsilon_s-Y^\delta_s|^2+E'|Y^{\epsilon\prime}_s-Y^{\delta\prime}_s|^2
 \nonumber
\\&&+|Z^\epsilon_s-Z^\delta_s|^2+E'|Y^{\epsilon\prime}_s-Y^{\delta\prime}_s|^2]ds.
\end{eqnarray}
Since
\begin{eqnarray*}
\langle Y^\epsilon_s-Y^\delta_s,
\nabla\varphi_{\epsilon}(Y^\epsilon_s)-\nabla\varphi_{\delta}(Y^\delta_s)\rangle\geq
-(\epsilon+\delta)|\nabla\varphi_{\epsilon}(Y^\epsilon_s)||\nabla\varphi_{\delta}(Y^\delta_s)|,
\end{eqnarray*}
combining \eqref{lm1:18} and \eqref{lm1:19}, by the same procedure as the proof of
Lemma \ref{lemma:4.1}, there exists a constant $C>0$ such that
\begin{eqnarray}\label{lm1:20}
E\int_0^Te^{\lambda
s}(|Y^\epsilon_s-Y^\delta_s|^2+|Z^\epsilon_s-Z^\delta_s|^2)ds\leq
C(\epsilon+\delta)M_2,
\end{eqnarray}
and consequently, we can conclude the proof by Burkholder--Davis--Gundy's inequality and \eqref{lm1:20}.

Now, we aim to give the proof of Theorem \ref{theorem4:1}.
\\
\textbf{Proof.} {\bf Existence.}  From Lemma \ref{lemma:4.3}, we can deduece that there exist stochastic
  processes $Y\in S^{2}_{\mathbb{F}}(0,T;R)$ and $Z\in
H^2_{\mathbb{F}}(0,T;R^{d})$ such that
$$\lim_{\epsilon\rightarrow0}(Y^\epsilon, Z^\epsilon)=(Y, Z).$$
Consequently, Lemma \ref{lemma:4.2} implies that
\begin{eqnarray}
\lim_{\epsilon\rightarrow0}J_{\epsilon}(Y^\epsilon)=Y \ \ {\rm in}\
\ H^2_{\mathbb{F}}(0,T;R)\nonumber
\end{eqnarray}
and
\begin{eqnarray}
\lim_{\epsilon\rightarrow0}E[e^{\beta
t}|J_{\epsilon}(Y^\epsilon_t)-Y_t|^2]=0,\ \ 0\leq t\leq
T.\nonumber
\end{eqnarray}
Moreover, Fatou's lemma, (ii) of Lemma \ref{lemma:4.2}, Proposition
\ref{prop4:1} and the lower semicontinuity of $\varphi$ shows that
(ii) of Definition \ref{def:1} is satisfied.

On the other hand, (i) of Lemma \ref{lemma:4.2} shows that
$U^\epsilon_t:=\nabla\varphi_{\epsilon}(Y^\epsilon_t)$ are bounded
in the space $H^2_{\mathbb{F}}(0,T;R)$, so there exists a
subsequence $\epsilon_n\rightarrow 0$ such that
\begin{eqnarray*}
U^{\epsilon_n} \rightarrow U, \ \ {\rm weakly \ in}\ \
H^2_{T}(0,T;R).
\end{eqnarray*}
Furthermore, we have
\begin{eqnarray*}
E\int_0^T|U_s|^2ds\leq\liminf_{n\rightarrow\infty}E\int_0^T|U^{\epsilon_n}_s|^2ds\leq
CM_2.
\end{eqnarray*}
In virtue of (H1), by passing limit in mean-filed BSDE \eqref{mfbsde:3}, we
deduce that the triple $(Y,Z,U)$ satisfies (iv) of Definition
\ref{def:1}.

Finally, since
 $U^\epsilon_t\in\partial\varphi(J_{\epsilon}(Y^{\epsilon}_t))$,
 $t\in[0,T]$, it follows that, for all $V\in
 H^2_{\mathbb{F}}(0,T;R)$,
\begin{eqnarray}
e^{\beta t}\langle U^\epsilon_t,
V_t-J_{\epsilon}(Y^{\epsilon}_t)\rangle+e^{\beta
t}\varphi(J_{\epsilon}(Y^{\epsilon}_t))\leq e^{\beta
t}\varphi(V_t),\  dP\times dt-{\rm a.e.}\nonumber
\end{eqnarray}
Taking the $\liminf$ in the probability in the above inequality, then
(iii) of  Definition \ref{def:1} holds.
\\
{\bf Uniqueness.}    Let $(Y^i_t, Z^i_t, U^i_t)$, $i=1, 2$ be two
solutions of mean-filed BSDE \eqref{mfbsde:2}. We denote by
\begin{eqnarray}
(\Delta Y_t, \Delta Z_t, \Delta U_t):=(Y^1_t-Y^2_t,
Z^1_t-Z^2_t, U^1_t-U^2_t). \nonumber
\end{eqnarray}
Applying It\^{o}'s formula to $e^{\beta t}|\Delta Y_t|^2$ yields that
\begin{eqnarray}
&& e^{\beta t}|\Delta Y_t|^2+\int_t^Te^{\beta s}(\beta|\Delta
Y_s|^2+|\Delta Z_t|^2)ds
+2\int_t^Te^{\beta s}\langle\Delta Y_s,\Delta U_s\rangle ds
 \nonumber\\&=&
2\int_t^Te^{\beta s}\langle\Delta Y_s, E'[f(s,Y^{1\prime}_s,Z^{1\prime}_s,Y^1_s,
Z^1_s)-f(s,Y^{2\prime}_s,Z^{2\prime}_s,Y^2_s,
Z^2_s)]\rangle ds  \nonumber\\&&
-2\int_t^Te^{\beta s}\langle \Delta Y_s,\Delta Z_sdW_s\rangle.
\nonumber
\end{eqnarray}
Since
\begin{eqnarray*}
 \langle \Delta Y_s, \Delta U_s\rangle\geq 0,\ \
dP\times dt-{\rm a.e.}
\end{eqnarray*}
Thus, as the same procedure as the proof of Lemma \ref{lemma:4.1}, we can derive
the uniqueness of the solution. The proof is complete.

\section{Viscosity solution of a nonlocal parabolic variation inequality}

In this part, we will give a probability interpretation
 for the viscosity solutions of nonlocal parabolic
 variational inequalities via  mean-field BSDEs with
 subdifferential operator studied before.

Let us consider the following McKean-Vlasov SDE parameterized by the initial condition
 $(t,\zeta)\in [0,T]\times L^2(\Omega, \mathcal{F}_t, P; R^n)$:
\begin{eqnarray}\label{sde:2}
 \left\{ \begin{array}{l@{ }r}dX_s^{t,\zeta}=
E'[b(s, (X_s^{0,x_0})', X_s^{t,\zeta})]ds+E'[\sigma(s, (X_s^{0,x_0})', X_s^{t,\zeta})]dW_s, s\in[t,T],\\
X_t^{t,\zeta}=\zeta,
\end{array}\right.
\end{eqnarray}
where $b: \bar{\Omega}\times[0,T]\times R^n \times R^n\rightarrow R^n$ and
$ \sigma: \bar{\Omega}\times[0,T]\times R^n \times R^n\rightarrow R^{n\times d}$ are
two measurable functions satisfying the following assumptions:
\begin{enumerate}
  \item [\rm(H4)] $b(\cdot, 0, 0)$ and $\sigma(\cdot, 0, 0)$ are $\bar{\mathbb{F}}$-progressively measurable
  continuous processes and there exists some constant $C>0$ such that
  $$|b(t,x',x)|+|\sigma(t,x',x)|\leq C(1+|x|), {\rm a.s.,\ for\ all}\ 0\leq t\leq T, x,x'\in R^n;$$
  \item [\rm(H5)]$b$ and $\sigma$ are Lipschitz in $x,x'$, i.e., there is a constant $C>0$ such that

 $|b(t,x'_1,x_1)-b(t,x'_2,x_2)|+|\sigma(t,x'_1,x_1)-\sigma(t,x'_2,x_2)|$
 $$\leq C(|x'_1-x'_2|+|x_1-x_2|),
  {\rm a.s., \ for\ all}\ 0\leq t\leq T, x_1,x'_1,x_2,x'_2\in R^n.$$
\end{enumerate}
It is known that, under the assumptions (H4) and (H5), SDE \eqref{sde:2} has a unique strong
solution (see, e.g, \cite{Buckdahn2}). Moreover, it holds true that, for any $p\geq2$, there exists
$C_p\in R$ such that, for any $t\in[0,T]$ and $\zeta, \zeta'\in L^p(\Omega, \mathcal{F}_t, P, R^n)$,
\begin{eqnarray}
  \begin{array}{l@{ }r}\displaystyle E\left[\sup_{t\leq s\leq T}|X_s^{t,\zeta}-X_s^{t,\zeta'}|^p|\mathcal{F}_t\right]
  \leq C_p|\zeta-\zeta'|^p, \ P-{\rm a.s.},
 \\ \\
\displaystyle E\left[\sup_{t\leq s\leq T}|X_s^{t,\zeta}|^p|\mathcal{F}_t\right]
  \leq C_p(1+|\zeta|^p), \ P-{\rm a.s.}
\end{array}
\end{eqnarray}
Let us give two real-valued functions $f(t,x',x,y',y,z)$ and $h(x',x)$,  which are assumed to satisfy
the following conditions (H6).
\begin{enumerate}
  \item [\rm(i)] $h:\bar{\Omega} \times R^n \times R^n\rightarrow R$ is an
  $\bar{\mathcal{F}}_t\otimes \mathcal {B}(R^n)$-measurable random variable and
  $f: \bar{\Omega}\times[0,T] \times R^n \times R^n\times R\times R\times R^d\rightarrow R$
  is a measurable process such that $f(\cdot,x',x,y',y,z)$ is $\bar{\mathbb{F}}$-adapted, for all
  $(x',x,y',y,z)\in R^n \times R^n\times R\times R\times R^d$.
  \item [\rm(ii)]  There exists a constant $C>0$ such that

 $|f(t,x'_1,x_1,y'_1,y_1,z_1)-f(t,x'_2,x_2,y'_2,y_2,z_2)|+|h(x'_1,x_1)-h(x'_2,x_2)|$
 $$\leq C(|x'_1-x'_2|+|x_1-x_2|+|y'_1-y'_2|+|y_1-y_2|+|y_1-y_2|),\
  {\rm a.s.},$$
   for all $0\leq t\leq T, x_1,x'_1,x_2,x'_2\in R^n, y_1,y'_1,y_2,y'_2\in R $ and $z_1, z_2\in R^d$.
   \item [\rm(iii)] $f$ and $h$ satisfy a linear growth condition, i.e., there exists some $C>0$ such that,
    for all $x',x\in R^n$,
   $$|f(t,x'x,0,0,0)|+|h(x',x)|\leq C(1+|x|+|x'|), {\rm a.s.}$$
   \item [\rm(iv)] $f(\bar{\omega},t,x',x,y',y,z)$ is continuous in $t$ for all $(x',x,y',y,z)$,
   $P(d\bar{\omega})$-a.s.
   \item [\rm(v)] $f(t,x',x,y',y,z)$ is nondecreasing with respect to $y'$.
   \item [\rm(vi)] There exists some $m\in \mathbb{N}^*$ such that
   $\varphi(E[h(X_T^{0,x_0},x)])\leq C(1+|x|^m)$, $\forall x\in R^n$.
\end{enumerate}
Now, we consider the following coupled mean-field BSDE with subdifferential operator:
\begin{eqnarray}\label{mfbsde:5}
 \left\{ \begin{array}{l@{ }r}-Y_s^{t,\zeta}ds+\partial\varphi(Y_s)ds
\ni
E'[f(s, (X_s^{0,x_0})', X_s^{t,\zeta}, (Y_s^{0,x_0})',Y_s^{t,\zeta},Z_s^{t,\zeta})]ds-
 Z_s^{t,\zeta}dW_s,\ s\in[t,T],\\
Y_T^{t,\zeta}=E'[ h((X_T^{0,x_0})', X_T^{t,\zeta})].
\end{array}\right.
\end{eqnarray}
We first consider the equation \eqref{mfbsde:5} for
$(t,\zeta)=(0,x_0)$: We know that there exists
 a unique solution $(Y^{0,x_0},Z^{0,x_0},U^{0,x_0})\in S^{2}_{\mathbb{F}}(0,T; R)\times
 H^2_{\mathbb{F}}(0,T; R^{d})
\times H^2_{\mathbb{F}}(0,T; R)$ to the mean-field BSDE \eqref{mfbsde:5}. Once we get
 $(Y^{0,x_0},Z^{0,x_0},U^{0,x_0})$, equation \eqref{mfbsde:5} become a classical BSDE with
 subdifferential operator
 whose coefficients $\tilde{f}(s, X_s^{t,\zeta},y,z)]
 =E'[f(s, (X_s^{0,x_0})', X_s^{t,\zeta}, (Y_s^{0,x_0})',y,z)]$  and
 $\tilde{h}(X_T^{t,\zeta})=E'[ h((X_T^{0,x_0})', X_T^{t,\zeta})]$
 satisfy the assumption (H0). Thus, from theorem \ref{lemma:2.1}, we know
that there exists a unique solution
$(Y^{t,\zeta},Z^{t,\zeta},U^{t,\zeta})\in S^{2}_{\mathbb{F}}(0,T; R)\times H^2_{\mathbb{F}}(0,T; R^{d})
\times H^2_{\mathbb{F}}(0,T; R)$ to the equation \eqref{mfbsde:5}.

Now, let us define the random field
\begin{eqnarray}\label{def:2}
u(t,x):=Y^{t,x}_s|_{s=t}, (t,x)
\in [0,T]\times R^n,
\end{eqnarray}
where $Y^{t,x}$ is the solution of mean-field BSDE \eqref{mfbsde:5} with $x\in R^n$ at the place
of $\zeta\in  L^2(\Omega, \mathcal{F}_t, P; R^n)$.

In this section, we aim to study the following nonlocal parabolic variation inequality (PVI in short):
\begin{eqnarray}\label{pvi:1}
 \left\{ \begin{array}{l@{ }r}\displaystyle\frac{\partial u(t,x)}{\partial
 t}+Au(t,x)
 +E[f(t, X_t^{0,x_0}, x, u(t,X_t^{0,x_0}), u(t,x),
 Du(t,x)\cdot E[\sigma(t,X_t^{0,x_0},x)])]\in \partial \varphi(u(t,x)),\\
\displaystyle u(T,x)=E[h(X_T^{0,x_0}, x)], x\in R^n,
\end{array}\right.
\end{eqnarray}
where
$$Au(t,x):=\frac{1}{2}tr(a D^2u(t,x))
+\langle b, Du(t,x)\rangle$$
with  $a:= E[\sigma(t, X_t^{0,x_0}, x)]E[\sigma(t,
X_t^{0,x_0}, x)]^T$, $b:=E[b(t, X_t^{0,x_0}, x)]$.
Here the functions $b, \sigma, f$ and $h$
are supposed to satisfy (H4), (H5) and (H6), respectively, and $X^{0,x_0}$
is the solution of the SDE \eqref{sde:2}.
 Below, we denote by $\mathcal
{S}(n)$  the set of $n\times n$ symmetric non-negative
matrices.
\begin{Definition}\rm
Let $u\in C_p([0,T]\times R^n)$ and $(t,x)\in[0,T]\times R^n$. We
denote by $\mathcal{P}^{1,2,+}u(t,x)$: the set of triples $(p,q,X)\in
R\times R^n\times \mathcal {S}(n)$ which are such that
\begin{eqnarray*}
u(s,y)&\leq& u(t,x)+p(s-t)+\langle q, y-x\rangle
\\&&+\frac{1}{2}\langle X(y-x), y-x\rangle+o(|s-t|+|y-x|^2).
\end{eqnarray*}
$\mathcal{P}^{1,2,-}u(t,x)$ is defined similarly as the set of triples
$(p,q,X)\in R\times R^n\times \mathcal {S}(n)$ which are such that
\begin{eqnarray*}
u(s,y)&\geq& u(t,x)+p(s-t)+\langle q, y-x\rangle
\\&&+\frac{1}{2}\langle X(y-x), y-x\rangle+o(|s-t|+|y-x|^2).
\end{eqnarray*}
\end{Definition}

\begin{Remark}\rm
Here $C_p([0,T]\times R^n)=\Big\{u\in C([0,T]\times R^n)$: There
exists some constant $p>0$ such that  $\sup_{(t,x)\in [0,T]\times
R^n}\frac{|u(t,x)|}{1+|x|^p}<+\infty$\Big\}.
\end{Remark}

Next, we want to prove that $u(t,x)$ introduced by \eqref{def:2} is
the unique viscosity solution of PVI \eqref{pvi:1}. Before this, we first introduce the
definition of viscosity solution of PVI \eqref{pvi:1}, one can see Crandall, Ishii and Lions \cite{CIL1992}
for more details.
\begin{Definition} \rm
A random field $u\in C_p([0,T]\times R^n)$ which satisfies
$u(T,x)=E[h(X_T^{0,x_0}, x)]$.
\begin{enumerate}
  \item [\rm(i)]  $u$ is a viscosity subsolution
of PVI \eqref{pvi:1} if $$u(t,x)\in {\rm Dom}\varphi,\ \forall (t,x)\in
[0,T]\times R^n$$ and at any point $(t,x)\in [0,T]\times R^n$, for
any $(p,q,X)\in \mathcal{P}^{1,2,+}u(t,x)$, it holds that
\begin{eqnarray*}
-p-\frac{1}{2}{tr}(aX)-\langle b, q\rangle
-E[f(t, X_t^{0,x_0},x,u(t, X_t^{0,x_0}),u(t,x), q\cdot E[\sigma(t, X_t^{0,x_0}, x)])]
 \leq -\varphi'_{-}(u(t,x)).
\end{eqnarray*}
  \item [\rm(ii)] $u$ is a viscosity supsolution of PVI \eqref{pvi:1}
if $$u(t,x)\in {\rm Dom}\varphi,\ \forall (t,x)\in [0,T]\times R^n$$ and
at any point $(t,x)\in [0,T]\times R^n$, for any $(p,q,X)\in
\mathcal{P}^{1,2,-}u(t,x)$, it holds that
\begin{eqnarray*}
-p-\frac{1}{2}{tr}(aX)-\langle b, q\rangle -E[f(t,
X_t^{0,x_0},x,u(t, X_t^{0,x_0}),u(t,x), q\cdot E[\sigma(t,
X_t^{0,x_0}, x)])]
 \geq -\varphi'_{+}(u(t,x)).
\end{eqnarray*}
  \item [\rm(iii)] $u$ is a viscosity solution of PVI \eqref{pvi:1} if it is both
a viscosity subsolution and a supersolution  of PVI \eqref{pvi:1}.
\end{enumerate}
\end{Definition}

We can now state the main results of this section.
\begin{Theorem}
Under the assumptions {\rm(H4)--(H6)}, the function $u(t,x)$ defined by \eqref{def:2}
 is the viscosity
solution of PVI \eqref{pvi:1}.
\end{Theorem}
\textbf{Proof.} For each $(t,x)\in[0,T]\times R^n$, and $\epsilon\in]0,1]$, let
$(Y_s^{t,x,\epsilon}, Z_s^{t,x,\epsilon}), s\in[t,T]$, the solution
of the mean-field BSDE
\begin{eqnarray}\label{mfbsde:6}
 \left\{ \begin{array}{l@{ }r}\displaystyle Y_s^{t,x,\epsilon}+\int_s^T\nabla\varphi_\epsilon(Y_s^{t,x,\epsilon})dr
=\int_s^T
E'[f(r, (X_r^{0,x_0})', X_r^{t,x}, (Y_r^{0,x_0,\epsilon})',Y_r^{t,x,\epsilon},Z_r^{t,x,\epsilon})]dr
-\int_s^TZ_r^{t,x,\epsilon}dW_r,\\
\displaystyle Y_T^{t,x,\epsilon}=E'[h((X_T^{0,x_0})', X_T^{t,x})].
\end{array}\right.
\end{eqnarray}
It is known from Buckdahn et al. \cite{Buckdahn2} that
$$u_{\epsilon}(t,x):=Y_t^{t,x,\epsilon},\ t\in[0,T],\ x\in R^n$$
is the unique viscosity solution (the authors, in \cite{Buckdahn2},
gave an example to explain why the uniqueness is only in $C_p([0, T]\times R^n)$) of the following
nonlocal PDE:
\begin{eqnarray}\label{pde:2}
 \left\{ \begin{array}{l@{ }r}\displaystyle \frac{\partial u_\epsilon(t,x)}{\partial
 t}+Au_\epsilon(t,x)
 +E[f(t, X_t^{0,x_0}, x, u_\epsilon(t,X_t^{0,x_0}), u_\epsilon(t,x),
 Du_\epsilon(t,x)\cdot E[\sigma(t,X_t^{0,x_0},x)])]= \nabla\varphi_\epsilon(u_\epsilon(t,x)),\\
\displaystyle u_\epsilon(T,x)=E[h(X_T^{0,x_0}, x)], x\in R^n.
\end{array}\right.
\end{eqnarray}
Moreover, it follows from Lemma \ref{lemma:4.3} that
\begin{eqnarray}\label{conv:1}
|u_\epsilon(t,x)-u(t,x)|\rightarrow 0, \  {\rm as}\   \epsilon
\rightarrow 0
\end{eqnarray}
for all $(t,x)\in[0,T]\times R^n$.

Let's first show that $u$ is the subsolution of PVI \eqref{pvi:1}. For
$(t,x)\in[0,T]\times R^n$ and $(p,q,X)\in \mathcal{P}^{1,2,+}u(t,x)$,
it follows from Crandall, Ishii and Lions \cite{CIL1992} that there
exist sequence
\begin{eqnarray}
 \left\{ \begin{array}{l@{ }r}\epsilon_n\searrow 0,\\(t_n,x_n)\in[0,T]\times R^n,\\
(p_n,q_n,X_n)\in \mathcal{P}^{1,2,+}u_{\epsilon_n}(t_n,x_n)
\end{array}\right.
\end{eqnarray}
such that
$$(t_n,x_n,u_{\epsilon_n}(t_n,x_n),p_n,q_n,X_n)\rightarrow (t,x,u(t,x),p,q,X) \ \ {\rm as} \ \
 n\rightarrow+\infty.$$
But for any $n$, we have
\begin{eqnarray}\label{ineq:61}
&&-p_n-\frac{1}{2}{tr}(a_nX_{n})-\langle b_n, q_n\rangle -E\left[f\left(t_n,
X_{t_n}^{0,x_0},x_n,u_{\epsilon_n}(t_n,
X_{t_n}^{0,x_0}),u_{\epsilon_n}(t_n,x_n), q_n\cdot E\left[\sigma(t_n,
X_{t_n}^{0,x_0}, x_n)\right]\right)\right]
 \nonumber\\&&\leq -\nabla\varphi_{\epsilon_n}(u_{\epsilon_n}(t_n,x_n)).
\end{eqnarray}
where $a_n:=E[\sigma(t_n, X_{t_n}^{0,x_0},
x_n)]E[\sigma(t, X_{t_n}^{0,x_0}, x_n)]^T$, $b_n:=E[b(t_n, X_{t_n}^{0,x_0}, x_n)]$.
Arguing as in Pardoux-R\c{a}\u{s}canu \cite{PardouxRascanu98}, we
let $y\in {\rm Dom}\varphi$ such that $y<u(t,x)$. Then by
\eqref{conv:1}, the uniformly convergence $u_\epsilon\rightarrow u$
on compacts implies that there exists $n_0>0$ such that
$y<u_{\epsilon_n}(t_n,x_n)$, $\forall n\geq n_0$. Thus, multiplying
both sides of \eqref{ineq:61} by $u_{\epsilon_n}(t_n,x_n)-y$, we get
\begin{eqnarray}\label{ineq:62}
&&\left\{-p_n-\frac{1}{2}{ tr}(a_nX_{n})
-\langle b_n, q_n\rangle
-E[f(t_n, X_{t_n}^{0,x_0},x_n,u_{\epsilon_n}(t_n,
X_{t_n}^{0,x_0}),u_{\epsilon_n}(t_n,x_n), q_n\cdot E[\sigma(t_n,
X_{t_n}^{0,x_0}, x_n)])]\right\}
\nonumber\\&&(u_{\epsilon_n}(t_n,x_n)-y)\leq
\varphi(y)-\varphi(J_{\epsilon_n}(u_{\epsilon_n}(t_n,x_n))).
\end{eqnarray}
Passing to $\liminf_{n\rightarrow+\infty}$ on both sides of \eqref{ineq:62}, we
have that for all $y<u(t,x)$,
\begin{eqnarray}\label{ineq:63}
&&\left\{-p-\frac{1}{2}{ tr}(aX)-\langle b, q\rangle-E[f(t, X_t^{0,x_0},x,u(t,
X_t^{0,x_0}),u(t,x), q\cdot E[\sigma(t, X_t^{0,x_0},
x)])]\right\}(u(t,x)-y) \nonumber\\&& \leq \varphi(y)-\varphi(u(t,x)),
\end{eqnarray}
it follows that
\begin{eqnarray}\label{ineq:64}
&&-p-\frac{1}{2}{ tr}(aX)-\langle b, q\rangle
-E[f(t, X_t^{0,x_0},x,u(t, X_t^{0,x_0}),u(t,x), q\cdot
E[\sigma(t, X_t^{0,x_0}, x)])] \nonumber\\&& \leq
-\varphi'_{-}(u(t,x)),
\end{eqnarray}
i.e., $u$ is a viscosity subsolution of PVI \eqref{pvi:1}.
 By the similar arguments we can show that $u$ is a viscosity supsolution of PVI \eqref{pvi:1},
 and thus, we complete the proof.

\begin{Theorem}
Under the assumptions {\rm(H4)--(H6)}, PVI \eqref{pvi:1} has a unique viscosity
solution.
\end{Theorem}
\textbf{Proof.} Here, we adopt the same arguments appeared in EI
Karoui et al. \cite{Karoui1997} and Pardoux and R\v{a}\c{s}canu
\cite{PardouxRascanu98}.

Suppose that $u$ is a subsolution and $v$ a supsolution of PVI \eqref{pvi:1} such that
$u(T,x)=v(T,x)=E[h(X_T^{0,x_0}, x)], x\in R^n$.

Define
\begin{eqnarray*}
\bar{u}(t,x)&:=&u(t,x)e^{\lambda t}\xi^{-1}(x),\\
\bar{v}(t,x)&:=&v(t,x)\ e^{\lambda t}\xi^{-1}(x)+\frac{\epsilon}{t},
\\ \eta(x)&:=& \xi^{-1}(x)D\xi(x)=p(1+|x|^2)^{-1}x,\\
\kappa(t,x)&:=&\xi^{-1}(x)D^2\xi(x)=p(1+|x|^2)^{-1}I-p(p-2)(1+|x|^2)^{-2}x\otimes x,
\end{eqnarray*}
where $\xi(x):=(1+x^2)^{\frac{p}{2}}$.
  Then, it is straightforward that
 $\bar{u}$ and $\bar{v}$ satisfy that (we write below $u$, $v$ instead of $\bar{u}$, $\bar{v}$)
\begin{eqnarray}\label{pvi:2}
 &&-\frac{\partial u(t,x)}{\partial
 t}-\widetilde{A}u(t,x)-E[\widetilde{f}(t, X_t^{0,x_0}, x, u(t,X_t^{0,x_0}), u(t,x),
 Du(t,x)\cdot E[\sigma(t,X_t^{0,x_0},x)])]\nonumber \\
&\leq&  -e^{\lambda t}\xi^{-1}(x)\varphi'_{-}(e^{-\lambda t}\xi(x)u(t,x))
\end{eqnarray}
and
\begin{eqnarray}\label{pvi:3}
 &&-\frac{\partial v(t,x)}{\partial
 t}-\widetilde{A}v(t,x)-E[\widetilde{f}(t, X_t^{0,x_0}, x, v(t,X_t^{0,x_0}), v(t,x),
 Dv(t,x)\cdot E[\sigma(t,X_t^{0,x_0},x)])] \nonumber \\
&\geq&  \frac{\epsilon}{t^2}-e^{\lambda t}\xi^{-1}(x)\varphi'_{+}\left(e^{-\lambda t}\xi(x)\left(v(t,x)
-\frac{\epsilon}{t}\right)\right),
\end{eqnarray}
where
$$\widetilde{A}\psi:=A\psi+\langle a\eta, D\psi \rangle+\left[\frac{1}{2}{ tr}(a\kappa)+\langle b,\eta \rangle-\lambda\right]\psi,$$
$$\widetilde{f}(t, X_t^{0,x_0}, x, u(t,X_t^{0,x_0}), u(t,x),
 Du(t,x)\cdot E[\sigma(t,X_t^{0,x_0},x)])$$
 $$:=e^{\lambda t}\xi^{-1}(x)f\Big(t, X_t^{0,x_0}, x, e^{-\lambda t}\xi(x)u(t,X_t^{0,x_0}),
 e^{-\lambda t}\xi(x)u(t,x),
 $$$$e^{-\lambda t}\xi(x) Du(t,x)\cdot E[\sigma(t,X_t^{0,x_0},x)])
 +e^{-\lambda t} u(t,x)D\xi(x)\cdot E[\sigma(t,X_t^{0,x_0},x)]\Big).$$
Let
$$F(t, X_t^{0,x_0}, x, u(t,X_t^{0,x_0}), u(t,x),
 Du(t,x), Du^2(t,x))$$
 $$:=-\widetilde{A}u(t,x)-E[\widetilde{f}(t, X_t^{0,x_0}, x, u(t,X_t^{0,x_0}), u(t,x),
 Du(t,x)\cdot E[\sigma(t,X_t^{0,x_0},x)])].$$
Then, by the Lipschitz condition (ii) of Assumption (H6), for $y_1>y_2$, we have
\begin{eqnarray*}
 &&\widetilde{f}(t,x',x, r, y_1,\theta)-\widetilde{f}(t,x', x, r, y_2,\theta) \nonumber \\
&=& e^{\lambda t}\xi^{-1}(x)\left[f(t, x', x, r, e^{-\lambda t}\xi(x)y_1,\theta)
 -f(t,x',x, r, e^{-\lambda t}\xi(x)y_2,\theta)\right]
\nonumber \\
&\leq&
C(y_1-y_2).
\end{eqnarray*}
Hence,
$$F(t,x',x, r, y_1,\mu, \nu)-F(t,x',x, r, y_2,\mu, \nu)
\geq \left[\lambda-\frac{1}{2}{ tr}(a\kappa)-\langle b,\eta
\rangle-C\right](y_1-y_2).$$ Since $a\kappa$ and $\langle b,\eta
\rangle$ are bounded, then we can choose $\lambda$ large enough such
that
$$y\rightarrow F(t, x', x, r, y,\mu, \nu)$$
is strictly increasing for any $(t, x', x, r, \mu, \nu)\in [0,T]
\times R^n \times R^n\times R\times R\times S(n)$, and thus $F$ is proper in the terminology of
\cite{CIL1992}.

What we need to show is that for any $K>0$, if $B_K:=\{|x|<K\}$,
$$\sup_{(0,T)\times B_K}(u-v)^+\leq \sup_{(0,T)\times {\partial B}_K}(u-v)^+.$$
Since the right-hand side tends to zero as $K\rightarrow \infty$, we will prove this fact by contradiction.

Assume that there exists some $K>0$ such that for some $(t_0, x_0)\in (0,T)\times B_K$
$$\delta:=u(t_0, x_0)-v(t_0, x_0)=\sup_{(0,T)\times B_K}(u-v)^+ > \sup_{(0,T)\times {\partial B}_K}(u-v)^+.$$

We define $(\hat{t},\hat{x},\hat{y})$ as being a point in
$[0,T]\times \bar{B}_K\times \bar{B}_K$ where the function
$$\Phi_\alpha(t,x,y)=u(t,x)-v(t,x)-\frac{\alpha}{2}|x-y|^2$$
achieves its maximum. Then by Lemma 8.7 in \cite{Karoui1997}, we have:

(i) for $\alpha$ large enough, $(\hat{t},\hat{x},\hat{y})\in [0,T]\times B_K\times B_K$,

(ii) $\alpha|\hat{x}-\hat{y}|^2\rightarrow 0$ and $|\hat{x}-\hat{y}|^2\rightarrow 0$,

(iii) $u(\hat{t},\hat{x})\geq v(\hat{t},\hat{y})+\delta$.\\
Then for $\alpha$ large enough,
$$e^{-\lambda \hat{t}}\xi(\hat{x})u(\hat{t},\hat{x})
\leq e^{-\lambda \hat{t}}\xi(\hat{y})\left(v(\hat{t},\hat{y})-\frac{\epsilon}{\hat{t}}\right)$$
and, as a result,
\begin{eqnarray}\label{ineq:4.1}
 -\varphi'_{-}(e^{-\lambda \hat{t}}\xi(\hat{x})u(\hat{t},\hat{x}))
\leq -\varphi'_{+}\left(e^{-\lambda \hat{t}}\xi(\hat{y})\left(v(\hat{t},\hat{y})-\frac{\epsilon}{\hat{t}}\right)\right).
\end{eqnarray}
Furthermore, from Theorem 8.3 in \cite{CIL1992}, we know that there exists
$$(p,X,Y)\in R\times \mathcal{S}(n)\times \mathcal{S}(n)$$
such that
$$(p,\alpha(\hat{x}-\hat{y}),X)\in  \mathcal{P}^{1,2,+}u(t,x),$$
$$(p,\alpha(\hat{x}-\hat{y}),Y)\in  \mathcal{P}^{1,2,-}v(t,x).$$
Next, because $u$ (resp. $v$) is a subsolution (resp. supsolution) and $F$ is proper, then
 following the same line as in \cite{Karoui1997}, it follows from  \eqref{ineq:4.1}
 that
$$F(\hat{t}, X_{\hat{t}}^{0,x_0}, \hat{y}, u(\hat{t},X_{\hat{t}}^{0,x_0}), v(\hat{t},\hat{y}),
 \alpha(\hat{x}-\hat{y}), Y)-F(\hat{t}, X_{\hat{t}}^{0,x_0}, \hat{x}, u(\hat{t},X_{\hat{t}}^{0,x_0}), v(\hat{t},\hat{y}),
 \alpha(\hat{x}-\hat{y}), X)
 \geq  \frac{\epsilon}{t^2}.$$
Finally, by the Lipschitz condition (ii) of (H6) on $f$, following the proof on Page 734 in \cite{Karoui1997},
 we can deduce a similar contradiction. The uniqueness is proved.


\end{document}